\documentclass{article}
\usepackage{amsmath, amsthm, amssymb, amscd, latexsym, url}

\begin{document}

\def\aec{{\mathcal K}}
\def\monst{{\mathfrak C}}
\def\cat{{\bf C}}
\def\dcat{{\bf D}}
\def\grp{{\bf Grp}}
\def\posetc{{\bf Pos}}
\def\comppos{{\bf CPO}}
\def\sets{{\bf Set}}
\def\topcat{{\bf Top}}
\def\rels{{\bf Rel(\Sigma)}}
\def\alphm{\alpha_\monst}
\def\autc{\mbox{Aut}(\monst)}
\def\ksst{{\prec}_{\aec}}
\def\csst{{\prec}_{\cat}}
\def\kpsst{{\prec}_{\aec'}}
\def\kspst{{\succ}_{\aec}}
\def\sst{\subseteq}
\def\lsst{\subseteq_L}
\def\spst{\supseteq}
\def\lamsub{{\mbox{Sub}_{\leq\lambda}(M)}}
\def\lamsuf{{\mbox{Sub}_{\leq\lambda}({f[M]})}}
\def\lamsup{{\mbox{Sub}_{\leq\lambda}({M'})}}
\def\musub{{\mbox{Sub}_{\leq\mu}(M)}}
\def\chisub{{\mbox{Sub}_{\leq\chi}(M)}}
\def\cktame{(\chi\mbox{,}\kappa)\mbox{-tame}}
\def\cltame{(\chi\mbox{,}\lambda)\mbox{-tame}}
\def\mltame{(\mu\mbox{,}\lambda)\mbox{-tame}}
\def\cleqktame{(\chi\mbox{,}\leq\kappa)\mbox{-tame}}
\def\cleqltame{(\chi\mbox{,}\leq\lambda)\mbox{-tame}}
\def\cinftame{(\chi\mbox{,}\infty)\mbox{-tame}}
\def\cinftamen{(\chi\mbox{,}\infty)\mbox{-tameness}}
\def\kltame{(\lambda\mbox{,}\kappa)\mbox{-tame}}
\def\kmtame{(\lambda\mbox{,}\kappa)\mbox{-tame}}
\def\klcomp{(\kappa,\lambda)\mbox{-compact}}
\def\oicomp{(\omega,\infty)\mbox{-compact}}
\def\iicomp{(\infty,\infty)\mbox{-compact}}
\def\kllocal{(\kappa,\lambda)\mbox{-local}}
\def\oppn{U_{p,N}}
\def\oppi{U_{p_i,N_i}}
\def\opf{U_{q'\restfn,f[N]}}
\def\opqnp{U_{q\restp,N'}}
\def\opqen{U_{q\restn,N}}
\def\topmla{X_M^{\lambda}}
\def\topfla{X_{f[M]}^{\lambda}}
\def\topnla{X_N^{\lambda}}
\def\toppla{X_{M'}^{\lambda}}
\def\topmm{X_M^{\mu}}
\def\topmk{X_M^{\kappa}}
\def\topma{X_{M_\alpha}^{\mu}}
\def\topmb{X_{M_\beta}^{\mu}}
\def\topmi{X_{M_i}^{\mu}}
\def\topmc{X_M^{\chi}}
\def\topmn{X_M^{\aleph_0}}
\def\unifml{U_M^{\lambda}}
\def\idmula{\mbox{Id}_{\mu,\lambda}}
\def\idmuch{\mbox{Id}_{\mu,\chi}}
\def\lsn{\text{LS(}{\aec}\text{)}}
\def\lsc{\text{LS(}{\cat}\text{)}}
\def\lsnp{\mbox{LS(}\aec'\mbox{)}}
\def\gatm{\mbox{ga-S(}M\mbox{)}}
\def\gatmi{\mbox{ga-S(}M_i\mbox{)}}
\def\gatma{\mbox{ga-S(}M_\alpha\mbox{)}}
\def\gatmo{\mbox{ga-S(}M_1\mbox{)}}
\def\gatmt{\mbox{ga-S(}M_2\mbox{)}}
\def\gatn{\mbox{ga-S(}N\mbox{)}}
\def\gatni{\mbox{ga-S(}N_i\mbox{)}}
\def\gatp{\mbox{ga-S(}M'\mbox{)}}
\def\gatmdp{\mbox{ga-S(}M''\mbox{)}}
\def\gatamn{\mbox{ga-tp(}a/M,N\mbox{)}}
\def\gatam{\mbox{ga-tp(}a/M\mbox{)}}
\def\gatbm{\mbox{ga-tp(}b/M\mbox{)}}
\def\tpam{\mbox{tp(}a/M\mbox{)}}
\def\tpbm{\mbox{tp(}b/M\mbox{)}}
\def\gatf{\mbox{ga-S(}f[M]\mbox{)}}
\def\gato{\mbox{ga-S(}N_1\mbox{)}}
\def\gatap{\mbox{ga-tp(}a/M'\mbox{)}}
\def\gatdp{\mbox{ga-tp(}a/M''\mbox{)}}
\def\gatdp{\mbox{ga-tp(}a/M''\mbox{)}}
\def\gatfan{\mbox{ga-tp(}f(a)/N\mbox{)}}
\def\gatann{\mbox{ga-tp(}a/N\mbox{)}}
\def\gatana{\mbox{ga-tp(}a/N_\alpha\mbox{)}}
\def\gatani{\mbox{ga-tp(}a_F/N_{\alpha_i}\mbox{)}}
\def\gatne{\mbox{ga-S(}N_\empt\mbox{)}}
\def\gatt{\mbox{ga-S(}N_2\mbox{)}}
\def\gatnsigi{\mbox{ga-S(}\nsigi\mbox{)}}
\def\gatnsig{\mbox{ga-S(}N_\sigma\mbox{)}}
\def\gattn{\mbox{ga-S(}\ntaun\mbox{)}}
\def\gattau{\mbox{ga-S(}N_\tau\mbox{)}}
\def\typem{\mbox{S(}M\mbox{)}}
\def\typea{\mbox{S(}A\mbox{)}}
\def\typeb{\mbox{S(}B\mbox{)}}
\def\ntaui{N_{\tau i}}
\def\ntauo{N_{\tau\restr 1}}
\def\ntaun{N_{\tau\restr n}}
\def\nsigi{N_{\sigma i}}
\def\nsigo{N_{\sigma\restr 1}}
\def\nsign{N_{\sigma\restr n}}
\def\psigi{p_{\sigma i}}
\def\psigij{p_{\sigma ij}}
\def\ptaunp{p_{\tau\restr (n+1)}}
\def\cplex{{\mathbb C}}
\def\autcf{\mbox{Aut}_{F}\mbox{(}\cplex\mbox{)}}
\def\autm{\mbox{Aut}_{M}\mbox{(}\monst\mbox{)}}
\def\autmi{\mbox{Aut}_{M_i}\mbox{(}\monst\mbox{)}}
\def\autn{\mbox{Aut}_{N}\mbox{(}\monst\mbox{)}}
\def\autsp{\mbox{Aut}_{N'}\mbox{(}\monst\mbox{)}}
\def\autp{\mbox{Aut}_{M'}\mbox{(}\monst\mbox{)}}
\def\autfm{\mbox{Aut}_{f[M]}\mbox{(}\monst\mbox{)}}
\def\autfn{\mbox{Aut}_{f[N]}\mbox{(}\monst\mbox{)}}
\def\rpm{\mbox{r}_{M',M}}
\def\rpf{\mbox{r}_{M',f[M]}}
\def\restn{\upharpoonright N}
\def\restm{\upharpoonright M}
\def\restmi{\upharpoonright M_i}
\def\restni{\upharpoonright N_i}
\def\resta{\upharpoonright M_\alpha}
\def\restp{\upharpoonright N'}
\def\restr{\upharpoonright}
\def\resto{\upharpoonright N_1}
\def\restt{\upharpoonright N_2}
\def\restfn{\upharpoonright f[N]}
\def\oppn{U_{p,N}}
\def\oppb{U_{p_\beta,N_\beta}}
\def\oppi{U_{p_i,N_i}}
\def\opnp{U_{p',N'}}
\def\opf{U_{q'|f[N],f[N]}}
\def\opqnp{U_{q|_{N'},N'}}
\def\opqen{U_{q\upharpoonright N,N}}
\def\pspl{p_{\lambda}\mbox{-space}}
\def\pspm{p_{\mu}\mbox{-space}}
\def\pspn{p_{\lsn}\mbox{-space}}
\def\ssq{\subseteq}
\def\spq{\supseteq}
\def\empt{\emptyset}
\def\invlim{\stackrel{\mbox{lim}}{\leftarrow}}
\def\fst{f^{\ast}}
\def\rml{\mbox{RM}^{\lambda}}
\def\rmm{\mbox{RM}^{\mu}}
\def\rmc{\mbox{RM}^{\chi}}
\def\rms{\mbox{RM}^{\lsn}}
\def\rll{\mbox{RM}^{\lambda}_{\ast}}
\def\rmo{\mbox{RM}^{\aleph_0}}
\def\rmn{\mbox{RM}^{\nu}}
\def\rmor{\mbox{RM}}
\def\cbl{\mbox{CB}^{\lambda}}
\def\kemb{\hookrightarrow_\aec}
\def\lsat{\lambda\mbox{-saturated}}
\def\lgsat{\lambda\mbox{-Galois-saturated}}
\def\lmh{\lambda\mbox{-model-homogeneous}}
\def\homc{\mbox{Hom}_\cat}
\def\homk{\mbox{Hom}_\aec}
\def\homkp{\mbox{Hom}_{\aec'}}
\def\homtk{\mbox{Hom}(-,K)}
\def\hom{\mbox{Hom}}
\def\lpres{\mbox{\bf Pres}_\lambda(\cat)}
\def\prsha{\mbox{\bf Set}^{{\mathcal A}^{{op}}}}
\def\prshm{\mbox{\bf Set}^{{M}^{{op}}}}
\def\prshkp{\mbox{\bf Set}^{{(\aec')}^{{op}}}}
\def\preslcat{\mbox{\bf Pres}_\lambda(\cat)}
\def\lplus{\lambda^{+}}
\def\elemt{{\bf Elem(}T{\bf )}}
\def\lang{{L}}
\def\sless{\vartriangleleft}
\def\slesseq{\trianglelefteq}
\def\mlset{(M,\lambda)\mbox{\bf -Set}}
\def\mlpset{(M,\lambda^+)\mbox{\bf -Set}}
\def\molpset{(M^{op},\lambda^+)\mbox{\bf -Set}}
\def\mset{M\mbox{\bf -Set}}
\def\moset{M^{op}\mbox{\bf -Set}}
\def\klang{{L}(\aec)}
\def\lkstruct{\klang{\bf{-Struct}}}
\def\lstruct{\lang{\bf{-Struct}}}
\def\sigstruct{L(T){\bf{-Struct}}}
\def\subcat{{\bf C}}
\def\basecat{{\bf B}}
\def\lwow{{L_{{\omega}_1,\omega}}}
\def\lkw{{L_{\kappa,\omega}}}
\def\lkl{{L_{\kappa,\lambda}}}
\def\linfw{{L_{\infty,\omega}}}
\def\lunc{{L(Q)}}
\def\lwowq{{L_{{\omega}_1,\omega}(Q)}}
\def\cof{\mbox{cf}}

\theoremstyle{plain}
\newtheorem{thm}{Theorem}[section]
\newtheorem{cor}[thm]{Corollary}
\newtheorem{lemma}[thm]{Lemma}
\newtheorem{prop}[thm]{Proposition}
\newtheorem{notat}[thm]{Notation}
\newtheorem{cla}[thm]{Claim}

\theoremstyle{definition}
\newtheorem{defn}[thm]{Definition}
\newtheorem{notation}[thm]{Notation}
\newtheorem{condition}[thm]{Condition}
\newtheorem{example}[thm]{Example}
\newtheorem{introduction}[thm]{Introduction}

\theoremstyle{remark}
\newtheorem{rmk}[thm]{Remark}

\include{header}

\numberwithin{thm}{section}     

\title{\Large {RANK FUNCTIONS AND PARTIAL STABILITY SPECTRA FOR TAME ABSTRACT ELEMENTARY CLASSES}}
\author{MICHAEL J. LIEBERMAN\footnote{The contents of this paper are drawn from the author's doctoral dissertation, completed at the University of Michigan under the supervision of Andreas Blass.  The author also acknowledges useful conversations with John Baldwin and Alexei Kolesnikov.}}

\maketitle

\begin{abstract}We introduce a family of rank functions and related notions of total transcendence for Galois types in abstract elementary classes.  We focus, in particular, on abstract elementary classes satisfying the condition know as tameness---currently regarded as a necessary condition for the development of a reasonable classification theory---where the connections between stability and total transcendence are most evident.  As a byproduct, we obtain a partial upward stability transfer result for tame abstract elementary classes stable in a cardinal $\lambda$ satisfying $\lambda^{\aleph_0}>\lambda$, a substantial generalization of a result of \cite{bkv}.\end{abstract}

\section{\large\textnormal{INTRODUCTION}}\label{intro}
For several decades, one of the overwhelming preoccupations of model theorists has been the pursuit of a classification theory for nonelementary classes, encompassing both the resolution of various categoricity conjectures generalizing Morley's theorem and the development of something resembling stability theory in the non-first-order context.  The impetus to study such questions is as much practical as theoretical: many natural mathematical objects defy first-order axiomatization (e.g. Banach spaces, Artinian commutative rings, and the complex numbers with exponentiation), and are thus beyond the scope of the results and methods of the classical theory.  On a theoretical level, the study of nonelementary classes offers an opportunity to obtain a deeper understanding of how classical stability theory actually works---which aspects are purely general, and which are inextricably linked to the unique structure of first-order logic.

The initial analysis of the model theory of more general logics (such as $\lkw$) was largely rooted in syntactic considerations, unfailingly identifying types with satisfiable sets of formulas, and relying on methods closely tailored to the peculiarities of the logics in question.  Abstract elementary classes (AECs) and their associated, fundamentally non-syntactic Galois types (both notions due to Shelah) represent a broad, unifying framework that supports a more uniform treatment of questions concerning categoricity and stability, transcending the focus on individual logics that limited the cope of earlier assays.  Essentially the category-theoretic hulls of elementary classes---we discard syntax entirely, but retain the essential diagrammatic properties of elementary embeddings---AECs are sufficiently general to encompass classes of models of sentences in $\linfw$, $\lunc$, and $\lwowq$, as well as homogeneous classes, but also seemingly rich enough in structure to permit the development of an interesting classification theory.

Naturally, with the added generality comes added difficulty: for AECs satisfying a variety of conditions intended to guarantee reasonable behavior, two decades of work have typically yielded only partial categoricity and stability transfer results.  In \cite{shejar} and \cite{sh600}, for example, Shelah has established a variety of broad structural results (including a version of his categoricity conjecture) for AECs with good or semi-good frames; that is, AECs equipped with forking relations reminiscent of those associated with stable or superstable theories in the classical realm.  In \cite{grvacat}, Grossberg and Vandieren obtain a large-scale categoricity transfer result for the considerably larger class of AECs in which the Galois types satisfy the important locality condition known as tameness using a toolkit that includes such classical standbys as indiscernibles and Morley sequences.  With these same methods, augmented with the notion of splitting, Grossberg and VanDieren also establish, in \cite{grva}, a partial stability spectrum result for tame AECs.  

Most recently, Baldwin, Kueker and VanDieren use a splitting argument to prove a striking upward stability transfer result for $\aleph_0$-Galois-stable AECs (Theorem 2.1 in \cite{bkv}): if $\kappa$ is a cardinal of uncountable cofinality and the AEC is Galois stable in every cardinal less than $\kappa$, it must be Galois stable in $\kappa$ as well.  We here introduce an entirely new set of vaguely classical methods: a family of rank functions for Galois types, which bear a certain resemblance to Morley rank, and use these to prove a generalization of the result of \cite{bkv}.  In particular, we show that this kind of transfer is can be triggered not merely by $\aleph_0$-Galois stability, but by stability in any cardinal $\lambda$ satisfying the inequality $\lambda^{\aleph_0}>\lambda$, and significantly weaken the assumption on stability below $\kappa$.

\section{\large\textnormal{PRELIMINARIES}}\label{secaecprelim}

We begin with a very brief introduction to AECs, Galois types, and a few relevant properties thereof.  Naturally this exposition will not be exhaustive; readers interested in further details may wish to consult \cite{baldwinbk} or \cite{grclassth}.  To begin:

\begin{defn}\label{defaec} Let $L$ be a finitary signature (one-sorted, for simplicity).  A class of $L$-structures equipped with a strong submodel relation, $(\aec,\ksst)$, is an {\it abstract elementary class (AEC)} if both $\aec$ and $\ksst$ are closed under isomorphism, and satisfy the following axioms:
\begin{description}
\item[A0] The relation $\ksst$ is a partial order.
\item[A1] For all $M$, $N$ in $\aec$, if $M\ksst N$, then $M\subseteq_{L}N$.
\item[A2] (Unions of Chains) Let $(M_\alpha | \alpha<\delta)$ be a continuous $\ksst$-increasing sequence.
\begin{enumerate}
\item $\bigcup_{\alpha<\delta} M_\alpha\in\aec$.
\item For all $\alpha<\delta$, $M_\alpha\ksst\bigcup_{\alpha<\delta} M_\alpha$.
\item If $M_\alpha\ksst M$ for all $\alpha<\delta$, then $\bigcup_{\alpha<\delta} M_\alpha\ksst M$.\end{enumerate}
\item[A3] (Coherence) If $M_0, M_1\ksst M$ in $\aec$, and $M_0\subseteq_L M_1$, then $M_0\ksst M_1$.
\item[A4] (Downward L\"owenheim-Skolem) There exists an infinite cardinal $\lsn$ with the property that for any $M\in\aec$ and subset $A$ of $M$, there exists $M_0\in\aec$ with $A\subseteq M_0\ksst M$ and $|M_0|\leq|A|+\lsn$.
\end{description}
\end{defn}

The prototypical example, of course, is the case in which $\aec$ is an elementary class---the class of models of a particular first-order theory $T$---and $\ksst$ is the elementary submodel relation, in which case $\lsn$ is, naturally, $\aleph_0+|L(T)|$.  

For any infinite cardinal $\lambda$, we denote by $\aec_\lambda$ the subclass of $\aec$ consisting of all models of cardinality $\lambda$ (with the obvious interpretations for such notations as $\aec_{\leq\lambda}$ and $\aec_{>\lambda}$).  We say that $\aec$ is $\lambda$-categorical if $\aec_\lambda$ contains only a single model up to isomorphism.  For $M, N\in\aec$, we say that a map $f:M\to N$ is a $\aec$-embedding (or, more often, a strong embedding) if $f$ is an injective homomorphism of $\klang$-structures, and $f[M]\ksst N$, that is; $f$ induces an isomorphism of $M$ onto a strong submodel of $N$.  In that case, we write $f:M\kemb N$.  We also adopt the following notational shorthand: for any $M\in\aec$ and any cardinal $\lambda$, we denote by $\lamsub$ the set of all strong submodels $N\ksst M$ with $|N|\leq\lambda$.

\begin{defn}\label{defapjep} Let $\aec$ be an AEC.
\begin{enumerate} 
\item We say that an AEC $\aec$ has the {\it joint embedding property (JEP)} if for any $M_1, M_2\in\aec$, there is an $M\in\aec$ that admits strong embeddings of both $M_1$ and $M_2$, $f_i:M_i\kemb M$ for $i=1, 2$.
\item We say that an AEC $\aec$ has the {\it amalgamation property (AP)} if for any $M_0\in\aec$ and strong embeddings $f_1:M_0\kemb M_1$ and $f_2:M_0\kemb M_2$, there are strong embeddings $g_1:M_1\kemb N$ and $g_2:M_2\kemb N$ such that $g_1\circ f_1=g_2\circ f_2$.
\end{enumerate}
\end{defn}

Notice that both properties hold in elementary classes, as a consequence of the compactness of first-order logic.  In this more general context, devised to subsume classes of models in logics without any compactness to fall back on, both appear as additional (and nontrivial) assumptions on the class.  We will work exclusively with AECs that satisfy both the joint embedding and the amalgamation properties.

It is not immediately clear what we might embrace as a suitable notion of type in AECs, given that we have dispensed with syntax, and removed ourselves to a world of abstract embeddings and diagrams thereof.  The best candidate---the Galois type---has its origins in the work of Shelah (first appearing in \cite{shclassth2} and \cite{sh394}).  Although Galois types can be defined in very general AECs, we restrict our attention to those with amalgamation and joint embedding.  In any AEC $\aec$ of this form, we may fix a monster model $\monst\in\aec$, and consider all models $M\in\aec$ as strong submodels of $\monst$ with $|M|<|\monst|$.  In this case, Galois types have a particularly simple characterization.

\begin{defn}\label{defgaltype} Let $M\in\aec$, and $a\in\monst$.  The {\it Galois type of $a$ over $M$}, denoted $\gatam$, is the orbit of $a$ in $\monst$ under $\autm$, the group of automorphisms of $\monst$ that fix $M$.  We denote by $\gatm$ the set of all Galois types over $M$.\end{defn}

In case $\aec$ is an elementary class with $\ksst$ as elementary submodel, the Galois types over $M$ correspond to the complete first-order types over $M$: 
$$\gatam=\gatbm\mbox{ if and only if }\tpam=\tpbm$$
In general, however, Galois types and syntactic types do not match up, even in cases when the logic underlying the AEC is clear (say, $\aec=\mbox{Mod}(\psi)$, with $\psi\in\lwow$).  We run through a series of basic definitions and notations:

\begin{defn}\label{galtypedefs}
\begin{enumerate}
\item We say that $\aec$ is {\it $\lambda$-Galois stable} if for every $M\in\aec_\lambda$, $|\gatm|=\lambda$.
\item For any $M$, $a\in\monst$, and $N\ksst M$, the {\it restriction of $\gatam$ to $N$}, which we denote by $\gatam\restn$, is the orbit of $a$ under $\autn$.  This notion is well-defined: the restriction depends only on $\gatam$, not on $a$ itself.
\item Let $N\ksst M$ and $p\in\gatn$.  We say that $M$ {\it realizes} $p$ if there is an element $a\in M$ such that $\gatam\restn=p$.  Equivalently, $M$ realizes $p$ if the orbit in $\monst$ corresponding to $p$ meets $M$.
\item We say that a model $M$ is {\it $\lambda$-Galois-saturated} if for every $N\ksst M$ with $|N|<\lambda$ and every $p\in\gatn$, $p$ is realized in $M$.
\end{enumerate}
\end{defn}

Henceforth, the word ``type'' should be understood to mean ``Galois type,'' unless otherwise indicated.  Moreover, when there is no risk of confusion, we will omit the word ``Galois'' altogether, speaking simply of types, $\lambda$-stability, and $\lambda$-saturation.

As in first-order classification theory, one of the central preoccupations of those working with AECs is the identification of ``dividing lines,'' in the sense of Shelah: properties of classes which, when present, guarantee nice structural properties, properties that disappear when one finds oneself on the wrong side of a dividing line, only to be replaced by new misbehaviors; that is, ``nonstructure.''  Certain of these properties for AECs (excellence, as described in \cite{grkol}, and the existence of good or semi-good frames, as considered in \cite{sh600} and \cite{shejar}, respectively) echo classical dividing lines (simplicity, stability, and superstability), but we here concern ourselves primarily with a property that has no particularly natural classical analogue: tameness of Galois types.  Roughly speaking, tameness means that types are completely determined by their restrictions to small submodels, a condition slightly reminiscent of the locality properties of syntactic types.  A possible intuition is the following: we may regard types over small models as playing a role analogous to formulas in first-order model theory, in which case tameness guarantees that types are determined by their constituent ``formulas."  In this analogy, tameness of Galois types may in fact be most closely identified with completeness in the realm of syntactic types.

We present two degrees of tameness:

\begin{defn}\label{defpairtame} 
\begin{enumerate} 
\item We say that $\aec$ is {\it $\chi$-tame} if for every $M\in\aec$, if $p$ and $p'$ are distinct types over $M$, then there is an $N\in\chisub$, such that $p\restn\neq p'\restn$.
\item We say that $\aec$ is {\it weakly $\chi$-tame} if for every saturated $M\in\aec$, if $p$ and $p'$ are distinct types over $M$, then there is an $N\in\chisub$, such that $p\restn\neq p'\restn$.
\end{enumerate}
\end{defn}

A few words about the significance of the assumption of tameness: although it holds automatically in elementary classes (as well as in homogeneous classes, among other settings), it fails in some contexts---\cite{bakol} and \cite{bash}, for example, construct examples of non-tame AECs from Abelian groups.  Moreover, it is independent from the other properties one associates with well-behaved AECs.  It is established in \cite{bash}, for example, that given an AEC with amalgamation, there is an AEC without amalgamation with precisely the same tameness spectrum.  Its utility more than justifies the loss of generality involved in its assumption, though: it plays an indispensable role in establishing all of the existing upward categoricity transfer results (such as those of \cite{grvaonesucc}, \cite{grvacat}, and \cite{sh394}), the downward transfer result included in Chapter 15 of \cite{baldwinbk}, and nearly all of the partial stability spectrum results of \cite{bkv}, \cite{grva}, and \cite{thesis}, although weak tameness occasionally suffices in \cite{bkv} and \cite{thesis}.  At present, attempts to produce results of this nature for non-tame classes have not met with any success.

We close with a brief topological aside (a more detailed treatment can be found in \cite{toppap}) which, although it will not necessarily be emphasized in the sequel, provides essential motivation for the definitions that follow.  Recall the identification that formed the basis of our intuition concerning tameness of Galois types---types over small submodels as ``formulas"---and fix our notion of ``small" as ``of cardinality less than or equal to $\lambda$," with $\lambda\geq\lsn$.  Just as the Stone space topology on syntactic types is generated from a basis of clopen sets, each one consisting of all complete types containing a given formula, we here consider sets of Galois types extending given ``formulas"; that is, for each model $M\in\aec$, we consider all sets of the form
$$\oppn=\{q\in\gatm\,:\,q\restn=p\}$$
where $N\in\lamsub$ and $p\in\gatn$.  It is an easy exercise to show that the $\oppn$ form a basis for a topology on $\gatm$---we denote the resulting space by $\topmla$---and, moreover, that they are clopen sets, analogous to the first-order case.  In fact, if $\aec$ is an elementary class (in which case Galois types may be identified with complete first-order types), we see that the topologies thus obtained refine the standard syntactic topology.

In fact, what we obtain is a spectrum of spaces $\topmla$ associated with each $M\in\aec$, with one space for each $\lambda\geq\lsn$, and a tight correspondence between topological properties of the spaces in the various spectra and model-theoretic properties of the AEC (see \cite{toppap} for a complete account).  In particular, $\chi$-tameness of an AEC means precisely that any two distinct types over a model $M\in\aec$ can be distinguished by their restrictions to a submodel of $M$ of size at most $\chi$; that is, they can be separated by basic open sets in each of the spaces $\topmla$ for $\lambda\geq\chi$.  It should come as no surprise, then, that we have:

\begin{prop}\label{tamespec} If $\aec$ is $\chi$-tame, $\topmla$ is Hausdorff for all $M\in\aec$ and $\lambda\geq\chi$.\end{prop}

Furthermore, in light of the fact that every $\topmla$ has a basis of clopens, the spaces in question will be totally disconnected---nearly Stone spaces.  As it happens, though, tameness rules out not merely compactness, but even countable compactness.  The chief complication results from the following result:
\begin{prop}\label{intprop} For any $M\in\aec$ and $\lambda\geq\lsn$, the intersection of any $\lambda$ open subsets of $\topmla$ is open.\end{prop}
\begin{proof} It suffices to show that the intersection of $\lambda$ basic open sets is open.  To that end, let $\{\oppi\,|\,i<\lambda\}$ be a family of open subsets of $\topmla$, and let $q\in\bigcap_{i<\lambda} \oppi$.  The union of the $N_i$ is of cardinality $\lambda$ so, by the Downward L\"owenheim Skolem Property, there is a submodel $N\kemb M$ that contains it.  By coherence, $N_i\kemb N$ for all $i$.  Consider the basic open neighborhood $\opqen$.  It certainly contains $q$.  For any $q'\in\opqen$, moreover, 
$$q'\restni=q'\restn\restni=q\restn\restni=q\restni=p_i$$
Hence $q'$ is contained in $\bigcap_{i<\lambda}\oppi$ and, as this holds for any such $q'$, $\opqen\sst\bigcap_{i<\lambda}\oppi$.\end{proof}

As a result, if an AEC is $\chi$-tame, the spaces $\topmla$ for $\lambda\geq\chi$ exhibit an extraordinary degree of separation, as captured in our final topological results, Propositions 6.5 and 7.1 in \cite{toppap}, respectively:

\begin{prop}\label{tamesep} Let $\aec$ be $\chi$-tame, and $\lambda\geq\chi$.  Any set $S\ssq\topmla$ with $|S|\leq\lambda$ is closed and discrete.\end{prop}

Thus there will be a number of infinite sets with no accumulation point, ruling out countable compactness.  We close with a closely-related characterization of what it means to be an accumulation point in $\topmla$:

\begin{prop}\label{tameacc} Let $\aec$ be $\chi$-tame, and $\lambda\geq\chi$.  A type $q\in\topmla$ is an accumulation point of a set $S\ssq\topmla$ if and only if for every neighborhood $U$ of $q$, $|U\cap S|>\lambda$.\end{prop}

\section{\large\textnormal{RANKS FOR GALOIS TYPES}}\label{rankdef}
We have just seen that if an AEC $\aec$ is $\chi$-tame, then for any $\lambda\geq\chi$ (and, of course, $\lambda\geq\lsn$), the criterion for a type $q\in\gatm$ to be an accumulation point of a family of types in $\topmla$ is quite stringent: every neighborhood must contain more than $\lambda$ types in the given family.  This condition inspires the definition of the following Morley-like ranks:
\begin{defn}\label{defrank}[$\rml$] Assume $\aec$ is $\chi$-tame for some $\chi\geq\lsn$.  For $\lambda\geq\chi$, we define $\rml$ by the following induction: for any $q\in\gatm$ with $|M|\leq\lambda$,
\begin{itemize}
\item $\rml[q]\geq 0$.
\item $\rml[q]\geq\alpha$ for limit $\alpha$ if $\rml[q]\geq\beta$ for all $\beta<\alpha$.
\item $\rml[q]\geq\alpha+1$ if there exists a structure $M'\kspst M$ such that $q$ has strictly more than $\lambda$ extensions to types $q'$ over $M'$ with the property that
$$\rml[q'\restn]\geq\alpha\mbox{ for all }N\in\lamsup$$
\end{itemize}
For types $q$ over $M$ of arbitrary size, we define
$$\rml[q]=\min\{\rml[q\restn]\,:\,N\in\lamsub\}.$$
\end{defn}

Before we proceed, a bit of motivation for this longwinded definition: in topologizing $\gatm$ as $\topmla$ (and in developing our intuition with regard to tameness), we essentially allow the types over substructures of size at most $\lambda$ to play a role analogous to that ordinarily played by formulas.  In defining the ranks $\rml$, we first define the rank of these ``formulas" and subsequently define the rank on types extending them as the minimum of the ranks of their constituent ``formulas,'' just as in the definition of Morley rank in classical model theory.

There is something to check here, though.  The types over models of cardinality at most $\lambda$ whose ranks were defined by the inductive clause are also covered by the second clause.  We must ensure that there is no possibility of conflict between the two.  Let $\rll$ denote the ranks assigned to such types using only the inductive clause.

\begin{lemma}\label{partialmonot} The ranks $\rll$ are monotone: for any $N, N'\in\aec_{\leq\lambda}$, $N'\ksst N$, and any $q\in\gatn$, $\rll[q]\leq\rll[q\restp]$.\end{lemma}
\noindent{\bf Proof:} We show that for any ordinal $\alpha$, $\rll[q]\geq\alpha$ implies $\rll[q\restp]\geq\alpha$.  We proceed by an easy induction on $\alpha$.  The zero case is trivial, by the first bullet point.  The limit cases follow from the induction hypothesis.  For the successor case, notice that if $\rll[q]\geq\alpha+1$, the extension $M\kspst N$ that witnesses this fact is also an extension of $N'$ and witnesses that $\rll[q\restp]\geq\alpha+1$.  \hspace{3 mm}$\Box$

So we needn't worry that the rank assigned to a type $q\in\gatm$ with $M\in\aec_{\leq\lambda}$ under the second clause ($\rml[q]=\min\{\rml[q\restn]\,:\,N\ksst M, |N|\leq\lambda\}$) will differ from the rank of $p$ under the first clause.  Hence the ranks $\rml$ are well-defined.  In fact, we can now restate the third bullet point above in a more compact form: for $q\in\gatm$ with $M\in\aec_{\leq\lambda}$,
\begin{itemize}
\item $\rml[q]\geq\alpha+1$ if there exists a structure $M'\kspst M$ such that $q$ has strictly more than $\lambda$ many extensions to types $q'$ over $M'$ with $\rml[q']\geq\alpha$.\end{itemize}
We will invariably use this formulation in the sequel.

We now consider the basic properties of the ranks $\rml$, beginning with an extension of our partial monotonicity result, Lemma~\ref{partialmonot}.
\begin{prop}\label{monot}[Monotonicity] If $M\ksst M'$ and $q\in\gatp$, then $\rml[q]\leq\rml[q\restm]$.\end{prop}
\noindent{\bf Proof:} Trivial, from the second clause of the definition.\hspace{3 mm}$\Box$

Bear in mind that ``$p$ is a restriction of $q$'' corresponds to the first-order ``$q\vdash p$,'' so this is indeed an appropriate analogue of the monotonicity result for the classical Morley rank.
\begin{prop}\label{inv}[Invariance] Let $f$ be an automorphism of the monster model $\monst$, let $q\in\gatm$, and let $M'=f[M]$.  Then the type $f[q]\in\gatp$ satisfies $\rml[f[q]]=\rml[q]$.\end{prop}
\noindent{\bf Proof:} We again proceed by induction, and once again the zero and limit cases are trivial.  It remains to show that whenever $\rml[q]\geq\alpha+1$, $\rml[f[q]]\geq\alpha+1$, with the converse following by replaying the argument with $f^{-1}$ in place of $f$.  As before, if $\rml[q]\geq\alpha+1$, then each $N\in\lamsub$ has an extension $M_N$ over which $q\restn$ has more than $\lambda$ many extensions of rank at least $\alpha$.  Notice that for any $N'\in\lamsup$, $N'=f[N]$ for some $N\in\lamsub$ and, moreover, that 
$$f[q]\restp=f[q]\restr f[N]=f[q\restn]$$
The extensions of $q\restn$ to $M_N$ of rank at least $\alpha$ map bijectively under the automorphism $f$ to extensions of $f[q]\restp$ to types over $f[M_N]$ and, by the induction hypothesis, these image types are also of rank at least $\alpha$.  As a result, for any $N'\in\lamsup$, $\rml[f[q]\restp]\geq\alpha+1$.  Thus $\rml[q]\geq\alpha+1$, as desired.\hspace{3 mm}$\Box$

Trivially,
\begin{prop}\label{locchar}[$\leq\lambda$-Local Character] For any $q\in\gatm$, there is an $N\in\lamsub$ with $\rml[q\restn]=\rml[q]$.\end{prop}

Slightly less trivially,
\begin{prop}\label{rmlrels}[Relations Between $\rml$] Whenever $\lambda\leq\mu$, $\rmm[q]\leq\rml[q]$.\end{prop}
\noindent{\bf Proof:} By induction, with the zero and limit cases still trivial.  If $\rmm[q]\geq\alpha+1$, then for every $N\in\musub$, $q\restn$ has more than $\mu$ extensions to types of $\rmm$-rank at least $\alpha$ over some $M_N\kspst N$.  In particular, there are more than $
\lambda$ such extensions and, by the induction hypothesis, they are all of $\rml$-rank at least $\alpha$.  Noticing that $\lamsub\ssq\musub$, one sees that the same holds for all $N\in\lamsub$, meaning that, ultimately, $\rml[q]\geq\alpha+1$.\hspace{3 mm}$\Box$

The ranks $\rml[-]$ also have the following attractive property: letting $\cbl[q]$ denote the topological Cantor-Bendixson rank of a type $q\in\topmla$, where $M$ is the domain of $q$, we have
\begin{prop}\label{boundcbl} For any $q$ (say $q\in\gatm$), $\cbl[q]\leq\rml[q]$.\end{prop}
\noindent{\bf Proof:} We show by induction that $\cbl[q]\geq\alpha$ implies $\rml[q]\geq\alpha$.  The zero and limit cases are not interesting.  If $\cbl[q]\geq\alpha+1$, $q$ is an accumulation point of types of $\cbl$-rank at least $\alpha$.  This means, by Proposition~\ref{tameacc}, that every basic open neighborhood of $q$ contains more than $\lambda$ such types, which is to say that each $q\restn$ with $N\in\lamsub$ has more than $\lambda$ many extensions to types over $M$ of $\cbl$-rank at least $\alpha$.  By the induction hypothesis, these types are of $\rml$-rank at least $\alpha$, and for each $N\in\lamsub$ witness the fact that $\rml[q\restn]\geq\alpha+1$.  Naturally, it follows that $\rml[q]\geq\alpha+1$.\hspace{3 mm}$\Box$

A quick fact about $\cbl$, incidentally, which we obtain in the usual way:
\begin{prop}\label{cblbddtoisodense} If every $q\in\topmla$ has ordinal $\cbl$-rank, isolated points are dense in $\topmla$.\end{prop}

In particular, then, if all of the types over a model $M\in\aec$ have ordinal $\rml$-rank, isolated points are dense in $\topmla$.

One of the great virtues of the classical Morley rank is that types have unique extensions of same Morley rank: if $p\in\typea$ has $\rmor[p]=\alpha$ and $B\spst A$, there is at most one type $q\in\typeb$ that extends $p$ and has $\rmor[q]=\alpha$.  While we do not have a unique extension property for the ranks $\rml$, we do have a very close approximation, Proposition~\ref{quext} below.  First, we notice:
\begin{lemma}\label{quextlemma} Let $M\ksst M'$, $q\in\gatm$ with an ordinal $\rml$-rank.  There are at most $\lambda$ extensions $q'\in\gatp$ with $\rml[q']=\rml[q]$.\end{lemma}
\noindent{\bf Proof:} For the sake of notational convenience, say $\rml[q]=\alpha$.  Let $S=\{q'\in\gatp\,|\,q'\restm=q,\mbox{ and }\rml[q']=\alpha\}$, and suppose that $|S|>\lambda$.  For any $N\in\lamsub$, each $q'$ is an extension of $q\restn$, meaning that $\rml[q\restn]\geq\alpha+1$ for all such $N$.  But then $\rml[q]\geq\alpha+1$.  Contradiction.\hspace{3 mm}$\Box$

It is worth noting here that we have made no use of tameness in the discussion above---not in the definition of the ranks $\rml$, nor in the proofs of the properties thereof---and have merely made use of the fact that $\lambda\geq\lsn$.  From now on, however, it will often be critically important that we have the ability to separate distinct types, leading us to restrict attention to those $\rml$ with $\lambda\geq\chi$, the tameness cardinal of the AEC at hand.  We will be very explicit in making this assumption whenever it is necessary, beginning with the following essential property of the $\rml$:
\begin{prop}\label{quext}[Quasi-unique Extension] Let $\aec$ be $\chi$-tame with $\chi\geq\lsn$, and let $\lambda\geq\chi$.  Let $M\ksst M'$, $q\in\gatm$, and say that $\rml[q]=\alpha$.  Given any rank $\alpha$ extension $q'$ of $q$ to a type over $M'$, there is an intermediate structure $M''$, $M\ksst M''\ksst M'$, $|M''|\leq |M|+\lambda$, and a rank $\alpha$ extension $p\in\gatmdp$ of $q$ with $q'\in\gatp$ as its unique rank $\alpha$ extension.\end{prop}
\noindent{\bf Proof:} Again, let $S=\{q''\in\gatp\,|\,q''\restm=q,\mbox{ and }\rml[q'']=\alpha\}$.  From the lemma, $|S|\leq\lambda$, meaning that $S$ is discrete (by tameness and Corollary 3.35).  Thus there is some $N\in\lamsup$ with the property that $q'\restn$ satisfies $U_{q'\restn,N}\cap S=\{q'\}$.  Take $M''\in\aec$ containing both $M$ and $N$ and with $|M''|\leq |M|+|N|+\lsn\leq|M|+\lambda$, and set $p=q'\restr M''$.  Notice that, by monotonicity, 
$$\alpha=\rml[q']\leq\rml[p]\leq\rml[q]=\alpha$$
so $\rml[p]=\alpha$.  Naturally, $q'$ is a rank $\alpha$ extension of $p$.  Any such extension $q''$ of $p$ is also a rank $\alpha$ extension of $q$, hence in $S$, and also an extension of $q'\restn$.  By choice of $N$, then, we must have $q''=q'$.\hspace{3 mm}$\Box$

A possible gloss for this complicated-looking result: while a type $p\in\gatm$ of $\rml$-rank $\alpha$ may have many extensions of rank $\alpha$ over a model $M'\kspst M$, we need only expand its domain ever so slightly (adding at most $\lambda$ elements of $M'$) to guarantee that the rank $\alpha$ extension is unique.  As innocuous as this proposition may seem, it is the linchpin of the analysis of stability in this framework, and will see extensive use in Sections~\ref{sectransf} and \ref{secwtame}.  Of course, quasi-unique extension applies only to types with ordinal $\rml$-rank.  To take the fullest possible advantage of this property, then, we would be wise to confine our attention to AECs where $\rml$ is ordinal-valued on all types associated with the class.

\section{\large\textnormal{STABILITY AND $\lambda$-TOTAL TRANSCENDENCE}}\label{secstabtt}

In classical model theory, we call a first-order theory totally transcendental if all types over subsets of the monster model are assigned ordinal values by the Morley rank.  We now define analogous notions for Galois types in AECs, one for each Morley-like rank $\rml$.
\begin{defn}\label{deftt} We say that $\aec$ is {\it $\lambda$-totally transcendental} if for every $M\in\aec$ and $q\in\gatm$, $\rml[q]$ is an ordinal.\end{defn}

\begin{prop}\label{ttup} If $\aec$ is $\lambda$-totally transcendental, it is $\mu$-totally transcendental for all $\mu\geq\lambda$.\end{prop}
\noindent{\bf Proof:} See Proposition~\ref{rmlrels}, the result concerning the relationship between the various Morley ranks.\hspace{3 mm}$\Box$

Putting this proposition together with Propositions~\ref{boundcbl} and~\ref{cblbddtoisodense}, we have:
\begin{prop}\label{ttisodense} If $\aec$ is $\lambda$-totally transcendental, then for any $M$ and $\mu\geq\lambda$, isolated points are dense in $\topmm$.\end{prop}

As a nice special case, if $\lsn=\aleph_0$, and $\aec$ is $\aleph_0$-tame and $\aleph_0$-totally transcendental, then isolated points are dense in each space $\topmn$.  This bears a certain resemblance to the classical result linking total transcendence to the density of isolated points in the spaces of syntactic types.

As interesting as this result may be, the true power of the notion of $\lambda$-total transcendence lies in the leverage it provides in rather a different area: bounding the number of types over models.  We will take advantage of this in Sections~\ref{sectransf} and \ref{secwtame} below to prove a variety of results related to Galois stability.  Before we give any further thought to $\lambda$-totally transcendental AECs, though, it seems natural to inquire whether the notion is, in fact, a meaningful one; that is, whether one can actually find $\lambda$-totally transcendental AECs.  One can, as the following proposition suggests:

\begin{thm}\label{stabtt} If $\aec$ is $\chi$-tame for some $\chi\geq\lsn$, and is $\lambda$-stable with $\lambda\geq\chi$ and $\lambda^{\aleph_0}>\lambda$, then $\aec$ is $\lambda$-totally transcendental.\end{thm}
\noindent{\bf Proof:} Suppose that $\rml$ is unbounded; that is, suppose that there is a type $p$ over some $N\ksst\monst$ with $\rml[p]=\infty$.  Indeed, we may assume that $|N|\leq\lambda$ (if there is a type of rank $\infty$ over a larger structure, consider one of its restrictions).  We proceed to construct a $\aec$-structure $\bar{N}$ of size $\lambda$ over which there are $\lambda^{\aleph_0}$ many types, thereby establishing that $\aec$ is not $\lambda$-stable.  We first need the following:
\begin{cla} For any $\lambda$, there exists an ordinal $\alphm$ such that for any $q$ over $M\ksst\monst$, $\rml[q]=\infty$ just in case $\rml[q]\geq\alpha$.\end{cla}
\noindent{\bf Proof:} The number of all such types can be bounded in terms of $|\monst|$, meaning that the set of their ranks cannot be cofinal in the class of all ordinals.\hspace{3 mm}$\Box{\mbox{ (Claim)}}$

{\noindent}We construct $\bar{N}$ as follows:

Step $1$: Set $p_{\empt}=p$.  Since $p_{\empt}$ satisfies $\rml[p_\empt]=\infty$, $\rml[p_\empt]>\alpha$.  Hence there is an extension $M\kspst N$ over which $p_\empt$ has more than $\lambda$ many extensions of rank at least $\alpha$ and thus of rank $\infty$.  Take any $\lambda$ many of them, say $\{q_i\in\gatm\,|\,i<\lambda\}$.  By Proposition~\ref{tamesep}, this set is discrete in $\topmla$, meaning that for each $i<\lambda$, there is an $N_i'\in\lamsub$ with the property that $q_j\restr N_i'\ne q_i\restr N_i'$ whenever $j\ne i$.  Let $N_\empt$ be a $\aec$-substructure of $M$ containing $N\cup(\bigcup_{i<\lambda}N_i)$, $|N_\empt|=\lambda$.  For each $i<\lambda$, define $p_i=q_i\restr N_\empt$.  Notice that:\\
$\bullet$ $\rml[p_i]\geq\rml[q_i]=\infty$\\
$\bullet$ We have $p_i\restr N=(q_i\restr N_\empt)\restr N=q_i\restr N=p_\empt$.\\
$\bullet$ If $p_i=p_j$, then 
$$q_i\restr N_i=(q_i\restr N_\empt)\restr N_i=p_i\restr N_i=p_j\restr N_i=(q_j\restr N_\empt)\restr N_i=q_j\restr N_i,$$ 
in which case we must have $i=j$.  That is, the $p_i$ are distinct. \\
So we have a family $\{p_i\in\gatne\,|\,i<\lambda\}$ of distinct extensions of $p_\empt$, all of rank $\infty$.

Step $n+1$: For each $\sigma\in\lambda^{n}$, we have a structure $N_\sigma$ of size at most $\lambda$ and a family of distinct types $\{\psigi\in\gatnsig\,|\,i<\lambda\}$, all of rank $\infty$.  For each $i<\lambda$, apply the process of step 1 to each $\psigi$.  By this method, we obtain (for each $i<\lambda$) an extension $\nsigi\kspst N_\sigma$ with $|\nsigi|=\lambda$ and a family $\{\psigij\in\gatnsigi\,|\,j<\lambda\}$ of distinct extensions of $\psigi$, all of rank $\infty$.

Step $\omega$: Notice that for each $\tau\in\lambda^\omega$, the sequence $N_\empt,\ntauo,\dots,\ntaun,\dots$ is increasing and, moreover, $\ksst$-increasing (by coherence and the fact that everything embeds strongly in $\monst$).  It follows that $N_\tau:=\bigcup_{n<\omega}\ntaun$ is in $\aec$.  Notice also that the $\ptaunp\in\gattn$ for $n<\omega$ form an increasing sequence of types over the structures in this union.  By $(\omega,\infty)$-compactness of AECs (see Theorem 11.1 in \cite{baldwinbk}), there is a type $p_\tau\in\gattau$ with the property that for all $n<\omega$, $p_\tau\restr\ntaun=p_{\tau\restr (n+1)}$.  Let $\bar{N}$ be a structure of size $\lambda$ containing the union
$$\bigcup_{\sigma\in\lambda^{<\omega}}N_\sigma$$
(which is possible, since the union is of size at most $\lambda\cdot\lambda^{<\omega}=\lambda\cdot\lambda=\lambda$).  For each $\tau\in\lambda^{\omega}$, let $q_\tau$ be an extension of $p_\tau$ to a type over $\bar{N}$.  There are $\lambda^{\aleph_0}$ many such types, all over $\bar{N}$---it remains only to show that they are distinct.  To that end, suppose that $\tau,\tau'\in\lambda^\omega$, $\tau\ne\tau'$.  It must be the case that for some $\sigma\in\lambda^{<\omega}$, $\tau=\sigma i\dots$ and $\tau'=\sigma j\dots$ with $i\ne j$.  Then
$$q_\tau\restr N_\sigma = (q_\tau\restr N_\tau)\restr N_\sigma = p_\tau\restr N_\sigma = p_{\sigma i}$$
Similarly, we have
$$q_{\tau'}\restr N_\sigma = (q_{\tau'}\restr N_{\tau'})\restr N_\sigma = p_{\tau'}\restr N_\sigma = p_{\sigma j}$$
Since $p_{\sigma i}$ and $p_{\sigma j}$ are distinct by construction, $q_\tau\ne q_{\tau'}$, and the types are distinct as claimed.  By our assumption that $\lambda^{\aleph_0}>\lambda$, then, $\aec$ is not stable in $\lambda$.\hspace{3 mm}$\Box$

So the notion of $\lambda$-total transcendence is far from vacuous: totally transcendental AECs exist.  Indeed, they exist in great abundance.

\section{\large\textnormal{TRANSFER RESULTS}}\label{sectransf}
As promised, we will now employ $\lambda$-total transcendence as a tool to analyze the proliferation of types over models in AECs.  In fact, we introduce two methods of bounding the number of types over large models, the first of which is captured by the following result:
\begin{thm}\label{dumbbound} Let $\aec$ be $\chi$-tame for some $\chi\geq\lsn$, and $\lambda$-totally transcendental with $\lambda\geq\chi+|\klang|$.  Then for any structure $M$ with $|M|>\lambda$, $|\gatm|\leq |M|^{\lambda}$.
\end{thm}
\noindent{\bf Proof:} Take $q\in\gatm$, $\rml[q]=\alpha$.  It must be the case that $\rml[q\restn]=\alpha$ for some $N\in\lamsub$ whence, by Proposition~\ref{quext}, there is an $N'\kspst N$ with $|N'|\leq\lambda$ such that $q\restr N'$ has rank $\alpha$ and has a unique rank $\alpha$ extension $q$ over $M$.  Hence
$$|\{q\in\gatm\,|\,\rml[q]=\alpha\}|\leq|\{p\mbox{ over }N\in\lamsub\,|\,\rml[p]=\alpha\}|$$
and thus
$$\begin{array}{rcl} |\gatm| & = & |\bigcup_\alpha \{q\in\gatm\,|\,\rml[q]=\alpha\}|\\
& \leq & |\bigcup_\alpha\{p\mbox{ over }N\in\lamsub\,|\,\rml[p]=\alpha\}|\\
& = & |\{p\mbox{ over }N\in\lamsub\}|\end{array}$$
Obviously, $|N|\leq\lambda$ for any $N\in\lamsub$, so by a general principle (see Proposition 2.10 in \cite{thesis}, for example) it must be the case that $|\gatn|\leq 2^\lambda$.  Hence
$$\begin{array}{rcl} |\gatm| & \leq & 2^\lambda\cdot |\lamsub|\\
& \leq & 2^\lambda\cdot |M|^\lambda\\
& = & |M|^\lambda\end{array}$$
So we are done.\hspace{3 mm}$\Box$

A major consequence of this theorem is:
\begin{cor}\label{dumbtransftt} If $\aec$ is $\chi$-tame for some $\chi\geq\lsn$, and $\lambda$-totally transcendental with $\lambda\geq\chi+|\klang|$, then for any $\mu$ satisfying $\mu^\lambda=\mu$, $\aec$ is $\mu$-stable.\end{cor}

In light of Theorem~\ref{stabtt}, we may eliminate any mention of $\lambda$-total transcendence and translate the corollary above into a pure upward stability transfer result:
\begin{prop}\label{dumbtransfstab} If $\aec$ is $\chi$-tame for some $\chi\geq\lsn$, and $\lambda$-stable with $\lambda\geq\chi+|\klang|$ and $\lambda^{\aleph_0}>\lambda$, then $\aec$ is stable in every $\mu$ with $\mu^{\lambda}=\mu$.\end{prop}
\noindent{\bf Proof:} Under the assumptions on $\lambda$, stability in $\lambda$ implies $\lambda$-total transcendence of $\aec$ (by the theorem invoked above).  We are then in the situation of Corollary~\ref{dumbtransftt}, and the result follows.\hspace{3 mm}$\Box$

It should be noted that this result is weaker than one found in \cite{grva}---namely Corollary 6.4---obtained by a splitting argument, which dispenses with the assumption that $\lambda^{\aleph_0}>\lambda$, inferring stability in $\mu$ with $\mu^\lambda=\mu$ from stability in any $\lambda$ larger than the tameness cardinal.  The potential advantages of our method become slightly more evident if we reinvent Proposition~\ref{dumbtransfstab} in such a way as to make clear the larger project for which it should be the jumping off point:
\begin{thm} Let $\aec$ be $\chi$-tame for some $\chi\geq\lsn$, and $\lambda$-stable with $\lambda\geq\chi+|L(\aec)|$ and $\lambda^{\aleph_0}>\lambda$.  Then there exists $\kappa\leq\lambda$ such that $\aec$ is $\mu$-stable in every $\mu$ satisfying $\mu^\kappa=\mu$.\end{thm}

Recall that if $\aec$ is stable in such a $\lambda$, it is $\lambda$-totally transcendental, and we may take
$$\kappa=\min\{\nu\,|\,\aec\mbox{ is }\nu\mbox{-totally transcendental}\}.$$
For this to be an actual improvement on existing results, of course, one would need to be able to give a meaningful bound on $\kappa$, the cardinal at which total transcendence first occurs, that places it well below $\lambda$.  Given that there is something of an analogy between this cardinal and the cardinal parameter of nonforking in elementary classes, there is some hope that this can be achieved, transforming the theorem above into a stability spectrum result reminiscent of those familiar from first-order model theory.

There is a different tack to be taken, though.  We introduce a second method for bounding the number of types over models of size larger than the cardinal of total transcendence:
\begin{thm}\label{bound} Let $\aec$ be $\chi$-tame for some $\chi\geq\lsn$, and $\lambda$-totally transcendental with $\lambda\geq\chi$.  For any $M\in\aec$ with $\mbox{cf(}|M|\mbox{)}>\lambda$, 
$$|\gatm|\leq|M|\cdot\sup\{|\gatn|\,|\,N\ksst M, |N|<|M|\}.$$\end{thm}
\noindent{\bf Proof:} Take a filtration of $M$:
$$M_0\ksst M_1\ksst\dots\ksst M_\alpha\ksst\dots\ksst M$$
for $\alpha<|M|$, with $|M_\alpha|<|M|$ for all $\alpha$ and $M=\bigcup_{\alpha<|M|} M_\alpha$.  Let $q\in\gatm$.  By $\lambda$-total transcendence, $\rml[q]=\beta$ for some ordinal $\beta$ and, moreover, this is witnessed by a restriction to a small submodel of $M$.  That is, there is a submodel $N\in\lamsub$ with $\rml[q\restn]=\beta$.  By Proposition~\ref{quext}, there is an intermediate extension $N\ksst N'\ksst M$ with $|N'|=\lambda$ such that $q\restp$ has unique rank $\beta$ extension $q$ over $M$.  Since $\mbox{cf(}|M|\mbox{)}>\lambda$, $N'\ssq M_\alpha$ for some $\alpha$ (and, by coherence, $N'\ksst M_\alpha$).  Clearly, the type $q\resta$ has $q$ as its only rank $\beta$ extension over $M$.  By a computation similar to that found in the proof of Proposition~\ref{dumbbound}, then, we have
$$|\gatm|\leq|M|\cdot\sup\{|\gatma|\,|\,\alpha<|M|\}$$
The inequality in the statement of the theorem is a trivial consequence.\hspace{3 mm}$\Box$

In essence, the theorem asserts that for models $M$ of cardinality $\mu$ with $\mbox{cf}(\mu)>\lambda$, the  equality $|\gatm|=\mu$ fails only if there is a submodel $N\ksst M$ with $|N|<\mu$ and $|\gatn|>\mu$.  To ensure that this does not occur---that we have, in short, stability in $\mu$---we need considerably less than full stability in the cardinals below $\mu$.  To be precise:
\begin{thm}\label{transftt} Let $\aec$ be $\chi$-tame for some $\chi\geq\lsn$ and $\lambda$-totally transcendental with $\lambda\geq\chi$, and let $\mu$ satisfy $\mbox{cf(}\mu\mbox{)}>\lambda$.  If for every $M\in\aec_{<\mu}$, $|\gatm|\leq\mu$, $\aec$ is $\mu$-stable.\end{thm}
\noindent{\bf Proof:} For any $M\in\aec_\mu$, 
$$|\gatm|\leq |M|\cdot\sup\{|\gatn|\,|\,N\ksst M, |N|<|M|\}\leq|M|\cdot |M|=|M|.\hspace{3 mm}\Box$$

We may tighten the statement up a bit, once we notice the following:
\begin{lemma}\label{transflemma} If there is a set $S$ of cardinals cofinal in an interval $[\kappa,\mu)$ which has the property that for every $M\in\aec$ with $|M|\in S$, $|\gatm|\leq\mu$, then $|\gatm|\leq\mu$ for every $M\in\aec_{<\mu}$.\end{lemma}
\noindent{\bf Proof:} Let $M\in\aec_{<\mu}$.  If $|M|\in S$, we are done.  If not, take $\nu\in S$ with $\nu>|M|$, and let $M'$ be a strong extension of cardinality $\nu$, $M\kemb M'\kemb\monst$.  By assumption, $|\gatp|\leq\mu$.  Every type over $M$ has an extension to a type over $M'$, and, of course, distinct types must have distinct extensions.  Hence 
$$|\gatm|\leq|\gatp|\leq\mu.\hspace{3 mm}\Box$$

So Theorem~\ref{transftt} becomes
\begin{thm}\label{besttransftt} Let $\aec$ be $\chi$-tame for some $\chi\geq\lsn$ and $\lambda$-totally transcendental with $\lambda\geq\chi$, and let $\mu$ satisfy $\mbox{cf(}\mu\mbox{)}>\lambda$.  If there is a set $S$ of cardinals cofinal in an interval $[\kappa,\mu)$ which has the property that for every $M\in\aec$ with $|M|\in S$, $|\gatm|\leq\mu$, $\aec$ is $\mu$-stable.\end{thm}

The assumption on the number of types over models of cardinality less than $\mu$ in the theorem above is very weak---it is certainly more than sufficient to assume that $\aec$ is $\nu$-stable for each $\nu$ in a cofinal set below $\mu$.  In particular, we have:
\begin{cor}\label{besttranfttwstab} If $\aec$ is $\chi$-tame for some $\chi\geq\lsn$ and $\lambda$-totally transcendental with $\lambda\geq\chi$, $\mu$ satisfies $\mbox{cf(}\mu\mbox{)}>\lambda$, and $\aec$ is stable on a set of cardinals cofinal in an interval $[\kappa,\mu)$, then $\aec$ is $\mu$-stable.\end{cor}

Using Theorem~\ref{stabtt} again, we may recast the corollary as a stability transfer result that makes no mention of total transcendence.  Though the product of this rewriting is a trivial consequence of the foregoing discussion (and, indeed, an extremely special case of Theorem~\ref{besttransftt}), it nonetheless generalizes a state-of-the-art result: Theorem 2.1 in \cite{bkv}.  First, the statement of the theorem:
\begin{thm}\label{besttransfstabwstab} If $\aec$ is $\chi$-tame for some $\chi\geq\lsn$ and $\lambda$-stable with $\lambda\geq\chi$ and $\lambda^{\aleph_0}>\lambda$, then for any $\mu$ with $\mbox{cf(}\mu\mbox{)}>\lambda$, if $\aec$ is stable on a set of cardinals cofinal in an interval $[\kappa,\mu)$, $\aec$ is $\mu$-stable.\end{thm}

If we restrict our attention to the special case in which the AEC is $\aleph_0$-stable, and assume more stability below $\mu$ than we need, strictly speaking, we have:
\begin{cor}\label{besttransfaleph0} If $\aec$ has $\lsn=\aleph_0$ and is $\aleph_0$-tame and $\aleph_0$-stable, then for any $\mu$ with $\mbox{cf(}\mu\mbox{)}>\aleph_0$, if $\aec$ is $\kappa$-stable for all $\kappa<\mu$, $\aec$ is $\mu$-stable.\end{cor}

In other words: in an AEC satisfying the hypotheses of the corollary, a run of stability will never come to an end at a cardinal $\mu$ of uncountable cofinality; it will, in fact, include $\mu$, then $\mu^+$, then $\mu^{++}$, and so on.  This remarkable fact is the aforementioned result of \cite{bkv} (which also appears as Theorem 11.11 in \cite{baldwinbk}).  What we have discerned by our method---and distilled in Theorem~\ref{besttransfstabwstab}---is that it is not stability in $\aleph_0$ which is critical, but rather stability in any cardinal $\lambda$ satisfying $\lambda^{\aleph_0}>\lambda$.  Moreover, it is not necessary to have stability in every cardinality below the cardinal of interest, $\mu$, but rather to have stability (or, indeed, somewhat less than stability) in a cofinal sequence of cardinals below $\mu$.

\section{\large\textnormal{WEAKLY TAME AECS}}\label{secwtame}

We now shift our frame of reference.  Thus far we have limited ourselves to AECs that are tame in some cardinal $\chi$.  The time has come to consider how the picture may differ in an AEC that is weakly $\chi$-tame, rather than $\chi$-tame.  The most important difference is that the argument for Theorem~\ref{stabtt} no longer works, meaning that we can no longer infer total transcendence from stability.  In a sense, then, we must treat total transcendence as a property in its own right, independent from the more conventional properties of AECs.  On the other hand, the argument for Theorem~\ref{bound} still goes through, provided the model in question is saturated.  That is,
\begin{prop}\label{wtamebd} Let $\aec$ be weakly $\chi$-tame for some $\chi\geq\lsn$, and $\mu$-totally transcendental with $\mu\geq\chi$.  If $M\in\aec$ is a saturated model with $\mbox{cf(}|M|\mbox{)}>\mu$, then $$|\gatm|\leq|M|\cdot\sup\{|\gatn|\,|\,N\ksst M, |N|<|M|\}$$\end{prop} 

We can use this to get bounds on the number of types over more general (that is, not necessarily saturated) models, provided that there are enough small saturated extensions in $\aec$.  In particular, if we can replace a model $M\in\aec_\lambda$ with a saturated extension $M'\in\aec_\lambda$, we have
$$|\gatm|\leq|\gatp|$$
where the latter is governed by the bound in the proposition above.  Given an adequate supply of such extensions, then, we can prove results similar to those found above, but now for weakly tame AECs.  The existence of saturated extensions also has the following consequence:
\begin{lemma}\label{exsatfewtypes} If every model in $\aec_\lambda$ has a saturated extension in $\aec_\lambda$, then for any $M\in\aec_{<\lambda}$, $|\gatm|\leq\lambda$.\end{lemma}
\noindent{\bf Proof:} Take $M\in\aec_{<\lambda}$.  Let $M'$ be a strong extension of $M$ of cardinality $\lambda$.  By assumption, $M'$ has a saturated extension $\bar{M}$ of cardinality $\lambda$.  This model, $\bar{M}$, realizes all types over submodels of size less than $\lambda$, hence realizes all types over the original model, $M$.  It follows that $|\gatm|\leq|\bar{M}|=\lambda$.\hspace{3 mm}$\Box$

Notice that this is precisely the type-counting condition from which we are able to infer stability using the bound in Proposition~\ref{wtamebd}.  Hence we have
\begin{thm}\label{wtametransf} Let $\aec$ be weakly $\chi$-tame for some $\chi\geq\lsn$, and $\mu$-totally transcendental with $\mu\geq\chi$.  Suppose that $\lambda$ is a cardinal with $\cof(\lambda)>\mu$, and that every $M\in\aec_\lambda$ has a saturated extension $M'\in\aec_\lambda$.  Then $\aec$ is $\lambda$-stable.\end{thm}
\noindent{\bf Proof:} Let $M\in\aec_\lambda$.  By assumption, we may replace $M$ with a saturated extension $M'$ of cardinality $\lambda$, over which there can be at most $\lambda$ types, by Proposition~\ref{wtamebd} and Lemma~\ref{exsatfewtypes}.  As noted above, $|\gatm|\leq|\gatp|$, and we are done.\hspace{3 mm}$\Box$

As a simple special case:
\begin{cor} Suppose $\aec$ is weakly $\chi$-tame for some $\chi\geq\lsn$ and $\mu$-totally transcendental with $\mu\geq\chi$, and suppose that $\lambda$ is a regular cardinal with $\lambda>\mu$.  If $\aec$ is stable on an interval $[\kappa,\lambda)$, $\aec$ is $\lambda$-stable.\end{cor}
\noindent{\bf Proof:} Let $M\in\aec$ be of cardinality $\lambda$.  By the usual union of chains argument, $M$ has a saturated extension $M'\in\aec$ which is also of cardinality $\lambda$.  Using the bound in Proposition~\ref{wtamebd}, we get $|\gatp|=\lambda$, whence $\lambda=|M|\leq|\gatm|\leq\lambda$, and thus $|\gatm|=\lambda$.  The result follows.\hspace{3 mm}$\Box$

In particular,
\begin{cor}\label{wtametransfsucc} Suppose $\aec$ is weakly $\chi$-tame for some $\chi\geq\lsn$, $\mu$-totally transcendental with $\mu\geq\chi$, and $\lambda$-stable for some $\lambda\geq\mu$.  Then $\aec$ is $\lambda^{+}$-stable.\end{cor}

This is weaker than a result of \cite{bkv}, which infers $\lambda^+$-stability from $\lambda$-stability in any AEC that is weakly tame in $\chi\leq\lambda$, with no additional assumptions.  The methods of \cite{bkv}---splitting, limit models---are better adapted to the task of proving local transfer results of this form in the context of weakly tame AECs.  The machinery of ranks and total transcendence provides us with leverage of a different sort, well suited to the task of producing partial spectrum results of a more global nature.  

In particular, we saw in Theorem~\ref{wtametransf} that in a $\mu$-totally transcendental AEC, even if merely weakly tame, we may infer stability in $\lambda$ given an adequate supply of saturated extensions of cardinality $\lambda$.  As it happens, the category-theoretic notion of weak $\lambda$-stability---introduced in \cite{rosdir}, and analyzed in the context of AECs in \cite{catpap}---ensures that this requirement is met.  The implications of this fact for the Galois stability spectra of weakly tame AECs are examined in \cite{catpap}.

\bibliographystyle{plain}
\bibliography{rankpaper2}
\end{document}